\documentclass{amsart}
\usepackage[all]{xy}
\usepackage{latexsym}
\usepackage{amssymb}
\usepackage{amsmath}
\usepackage{amsthm}
\usepackage{mathabx}
\usepackage{amscd}
\usepackage{fancyhdr}
\usepackage{fullpage}
\usepackage{mathrsfs}
\xyoption{arc}

 \usepackage{tikz}
\usetikzlibrary{matrix}
\usetikzlibrary{shapes}
\usetikzlibrary{arrows}
\usetikzlibrary{calc,3d}
\usetikzlibrary{decorations,decorations.pathmorphing,decorations.pathreplacing,decorations.markings}
\usetikzlibrary{through}

  \def\id{{\mbox{1 \hskip -7pt 1}}}
\newcommand{\sgn}{{\mathit s  \mathit g\mathit  n}}
 \newcommand{\lon}{\longrightarrow}
 \newcommand{\bu}{\bullet}
 
 \newcommand{\rar}{\rightarrow}

\newcommand{\Gra}{\mathcal{G} ra}

\newcommand{\fOGC}{\mathsf{fOGC}}

\newcommand{\Holie}{\mathcal{H} \mathit{olie}}
\newcommand{\Lie}{\mathcal{L} \mathit{ie}}

\newcommand{\Def}{\mathsf{Def}}
\newcommand{\fGC}{\mathsf{fGC}}

\newcommand{\OGC}{\mathsf{OGC}}
\newcommand{\fcOGC}{\OGC^0}

\newcommand{\hOGC}{\widehat{\OGC}}
\newcommand{\fhOGC}{\widehat{\fOGC}}
\newcommand{\GC}{\mathsf{GC}}

\newcommand{\fcGC}{\GC^0}

\newcommand{\fcGCdd}{\mathsf{{OGC}}_{d,d+1}}
\newcommand{\grt}{\fg\fr\ft}

 \newcommand{\bbu}{\mbox{\resizebox{2.8mm}{!}{$\bullet$}}}

\newcommand{\sGa}{{\mathsf{\Gamma}}}

 \newcommand{\Z}{{\mathbb Z}}
 \newcommand{\bS}{{\mathbb S}}

 \newcommand{\R}{{\mathbb R}}
 \newcommand{\N}{{\mathbb N}}
 \newcommand{\K}{{\mathbb K}}

 \newcommand{\ot}{\otimes}
 
 \newcommand{\Beq}{\begin{equation}}
 \newcommand{\Eeq}{\end{equation}}
 \newcommand{\Beqr}{\begin{eqnarray}}
 \newcommand{\Eeqr}{\end{eqnarray}}
 \newcommand{\Beqrn}{\begin{eqnarray*}}
 \newcommand{\Eeqrn}{\end{eqnarray*}}
 \newcommand{\Ba}{\begin{array}}
 \newcommand{\Ea}{\end{array}}
 \newcommand{\Bi}{\begin{itemize}}
 \newcommand{\Ei}{\end{itemize}}
 \newcommand{\Bc}{\begin{center}}
 \newcommand{\Ec}{\end{center}}

 \newcommand{\fg}{{\mathfrak g}}

\newcommand{\fr}{{\mathfrak r}}

\newcommand{\ft}{{\mathfrak t}}

\newcommand{\fC}{{\mathfrak C}}

 \newcommand{\cE}{{\mathcal E}}
 \newcommand{\cF}{{\mathcal F}}
 \newcommand{\cG}{{\mathcal G}}

 \newcommand{\ga}{\gamma}
 
 \newcommand{\Ga}{\Gamma}

 \newcommand{\sip}{\smallskip}
 \newcommand{\bip}{\bigskip}
 \newcommand{\mip}{\vspace{2.5mm}}

\theoremstyle{plain}
\swapnumbers

\newtheorem{prop-def}[theorem]{Proposition-definition}

\newtheorem{f-theorem}{Formality Theorem}[section]
\newtheorem{main-theorem}{Main~Theorem}[section]
\newtheorem{section-theorem}{Theorem}[section]

\theoremstyle{definition}

\begin{document}

 \sloppy

 \newenvironment{proo}{\begin{trivlist} \item{\sc {Proof.}}}
  {\hfill $\square$ \end{trivlist}}

\long\def\symbolfootnote[#1]#2{\begingroup%
\def\thefootnote{\fnsymbol{footnote}}\footnote[#1]{#2}\endgroup}

 \title{The oriented graph complex revisited}

\author{Sergei Merkulov}
\address{Sergei~Merkulov: Department of Mathematics, Luxembourg University,
Maison du Nombre, 6 Avenue de la Fonte,
 L-4364 Esch-sur-Alzette,   Grand Duchy of Luxembourg}
\email{sergei.merkulov@uni.lu}

\author{Thomas~Willwacher}
\address{Thomas~Willwacher: Department of Mathematics, ETH Zurich, Zurich, Switzerland}
\email{thomas.willwacher@math.ethz.ch}
\thanks{T.W. has been partially supported by the NCCR Swissmap, funded by the Swiss National Science Foundation.}

\author{Vincent Wolff}
\address{Vincent~Wolff: Department of Mathematics, Luxembourg University,
Maison du Nombre, 6 Avenue de la Fonte,
 L-4364 Esch-sur-Alzette,   Grand Duchy of Luxembourg}
\email{vincent.wolff@uni.lu}

 \begin{abstract} We prove that the Kontsevich graph complex $\GC_d^{2}$ and its oriented
  version $\OGC_{d+1}^2$  are  quasi-isomorphic as dg Lie algebras.


\end{abstract}
 \maketitle

{\large
\section{\bf Introduction}
}
\label{sec:introduction}

The famous graph complex $\GC_d$ was introduced in \cite{Ko1} in the context 
of deformation quantizations of Poisson structures. It is generated by connected graphs whose vertices are at least trivalent. The integer $d$ controls the cohomological degrees assigned to graphs, and the handling of their symmetries. 
The trivalence condition on Kontsevich's original graph complex may be relaxed. In the literature the larger complexes 
\[
\GC_d \subset \GC_d^{2}\subset \fcGC_d
\]
have been considered, with $\GC_d^{2}$ defined such that the vertices of the graphs are required to be at least bivalent, while for $\fcGC_d$ one drops the valence condition entirely.
All these complexes are in fact dg Lie algebras, the inclusions are compatible with the dg Lie structure, and the cohomological difference between them is controlled, see \S 2 below for a brief reminder.

The graph complexes above and their cohomology have seen many recent application in homological algebra and algebraic geometry. For example, they control the top weight cohomology of the moduli spaces of curves \cite{CGP1}, and also the homotopy automorphisms of the (rationalized) little $d$-disks operads \cite{FW}, and are hence of high interest to the community.

\sip

There is also an \emph{oriented} version $\fcOGC_d$ proposed by the first author. This is obtained by simply replacing the plain undirected graphs in Kontsevich's definition by directed acyclic graphs, that is, directed graphs with no closed paths of directed edges. Requiring vertices to be at least bivalent, one obtains the quasi-isomorphic dg Lie algebras
\[
\OGC_d^2\subset \OGC_d^0.
\]
Requiring in addition that there is at least one trivalent vertex one may define a yet smaller dg Lie subalgebra
\[
\OGC_d^3\subset \OGC_d^2.
\]
These oriented graph complexes have also been studied intensively in recent years. In particular, it has been proven in \cite{Wi2} that there is  an isomorphism  of Lie algebras at the cohomology level
 \begin{equation}\label{equ:fgc dgc}
 H^\bu(\fcGC_d)\simeq H^\bu(\fcOGC_{d+1})
 \end{equation}
for any $d\in \Z$. (Mind the shift in $d$.) Other simplified proofs of the above isomorphism --- but only at the 
level of graded vector spaces (i.e.\ not as Lie algebras) -- were given in \cite{Z1} and then in \cite{Me0}.
The oriented graph complexes also feature prominently in various mathematical problems: They appear in the description of the rational homotopy type of the space of long knots in $\R^n$ ($n\geq 4$) \cite{FTW}, they control the deformation theory of the PROP governing Lie bialgebras \cite{MWLieb}, and they appear centrally in the study of the real homotopy type of the real locus of the moduli space of genus zero curves \cite{KhoroshkinWillwacherMosaic}. 

\sip

The main result of this paper is an upgrade of the cohomology isomorphism \eqref{equ:fgc dgc} to a quasi-isomorphism on the level of dg Lie algebras. In the process we also provide yet another independent proof of the isomorphism \eqref{equ:fgc dgc}.

\sip

\subsection{Main Theorem}\label{1: Main Theorem} {\em For any $d\in \Z$ there exists a dg Lie 
algebra $\hOGC_{d,d+1}$ of graphs which fits into the following diagram of quasi-isomorphisms of 
dg Lie algebras}

$$
\fcGC_d    \stackrel{\pi_1}{\longleftarrow}   \hOGC_{d,d+1}  \stackrel{\pi_2}{\lon}
\fcOGC_{d+1}
$$
\emph{There is a grading by loop order on $\hOGC_{d,d+1}$ and the arrows $\pi_1$, $\pi_2$ respect the loop order gradings.}

\sip

The intertwining dg Lie algebra $\hOGC_{d,d+1}$ is explicitly described in \S 3: it is generated by graphs with two types of vertices and two types of edges, and is equipped
with a relatively non-trivial but explicitly described differential. This dg Lie algebra controls the deformation theory of a certain morphism
of 2-coloured operads as explained in \S 3 of this paper. The morphisms $\pi_1$ and $\pi_2$
are very simple and explicit. The theorem is proven in \S 4.

\sip

By extending the zigzag with the inclusions $\GC_d^2\subset \GC_d^0$ and $\OGC_d^2\subset \OGC_d^0$, which are known to be quasi-isomorphisms, we also obtain a zigzag of quasi-isomorphisms of dg Lie algebras
$$
\GC_d^{2} \xrightarrow{\sim} \bullet \xleftarrow{\sim} \OGC_{d+1}^2
$$
connecting the bivalent versions.
It is furthermore known that the inclusions $\GC_d\subset \GC_d^2$ and $\OGC_d^3\subset \OGC_d^2$ are quasi-isomorphisms in loop orders $\geq 2$. Hence extending our zigzag and truncating to loop orders $\geq 2$ we also obtain a zigzag of quasi-isomorphisms of dg Lie algebras
$$
\GC_d \xrightarrow{\sim} \bullet \xleftarrow{\sim} \OGC_{d+1}^3
$$

We note that we slightly deviate from the notation of the literature to obtain more consistent conventions for this paper. In particular the dg Lie algebra $\OGC^2_d$ has been denoted $\GC_d^{or}$ elsewhere.


\subsection{Some notation}   We work in this paper over a field $\K$ of characteristic zero.
 The set $\{1,2, \ldots, n\}$ is abbreviated to $[n]$;  its group of automorphisms is
denoted by $\bS_n$; the trivial (resp., the sign) one-dimensional representation of
 $\bS_n$ is denoted by $\id_n$ (resp.,  $\sgn_n$). The cardinality of a finite set $S$ is
 denoted by $\# S$ while its linear span over a
field $\K$ by $\K\left\langle S\right\rangle$.
The top degree skew-symmetric tensor power
of $\K\left\langle S\right\rangle$ is denoted by $\det S$; it is assumed that $\det S$
is a 1-dimensional Euclidean space associated with the unique Euclidean
structure on $\K\left\langle S\right\rangle$ in which the elements of $S$ serve
as an orthonormal basis; in particular, $\det S$ contains precisely two vectors of unit length.
If $V=\oplus_{i\in \Z} V^i$ is a graded vector space, then
$V[k]$ stands for the graded vector space with $V[k]^i:=V^{i+k}$. For $v\in V^i$ we set $|v|:=i$.


\bip

{\large
\section{\bf Recollections on graph complexes}
}
\label{2 sec:graph operad}

\sip

\subsection{An operad of graphs $\Gra_d$}\label{2: subsec on Gra_d}
Let $G_{n,l}$ be the set of graphs $\Ga$ with $n$ vertices (their set is denoted by $V(\Ga)$)
and $l$ directed edges (their set is denoted by $E(\Ga)$)  such that some bijections $V(\Ga)\rar [n]$ and $E(\Ga)\rar [l]$ are fixed, i.e.\ every edge and every vertex of $\Ga$ is marked.
Note that we allow tadpoles (self-edges) in our graph.
The permutation group (wreath product) $\bS_2\wr \bS_l$ acts on $G_{n,l}$ by relabeling the edges and flipping edge directions. For any
 any integer $d$ we define a collection of $\bS_n$-modules,
$$
\cG ra_{d}=\left\{\cG ra_d(n):= \prod_{l\geq 0} \K \langle G_{n,l}\rangle \ot_{ \bS_l}  \sgn_l^{\ot |d-1|} [l(d-1)] \sgn_2^{\otimes |d|l}  \right\}_{n\geq 1}
$$
where we assume that the group $\bS_n$ acts on $\Gra_d(n)$ by relabeling the vertices.

This is a $\Z$-graded vector space obtained by assigning to each edge of a generating graph $\Ga$ from $G_{n,l}$ the homological degree $1-d$ and providing $\Ga$ with a so called {\em orientation}\, $or$ which depends on the parity of $d$:  if $d$ is even, then $or$ a choice of ordering of edges up to an even permutation (an odd permutation acts as the multiplication by $-1$); if  $d$ is odd, then  $or$ is a choice
of directions on each edge (up to flip and the multiplication by $-1$).
\Beq\label{2: symmetry of dotted edges}
\Ba{c} \resizebox{11mm}{!}{\xy
 (10,1)*+{_k}*\frm{o}="B";
 (0,1)*+{_i}*\frm{o}="A";
 \ar @{.>} "A";"B" <0pt>
\endxy} \Ea
= (-1)^{d}
\Ba{c} \resizebox{11mm}{!}{\xy
 (10,1)*+{_k}*\frm{o}="B";
 (0,1)*+{_i}*\frm{o}="A";
 \ar @{<.} "A";"B" <0pt>
\endxy} \Ea.
\Eeq

\sip

This $\bS$-module
 is an operad with respect to the following operadic composition \cite{Ko2,Wi1},
$$
\Ba{rccc}
\circ_i: &  \cG ra_d(n) \times \cG ra_d(m) &\lon & \cG ra_d(m+n-1),
 \ \ \forall\ i\in [n]\\
         &       (\Ga_1, \Ga_2) &\lon &      \Ga_1\circ_i \Ga_2,
\Ea
$$
where  $\Ga_1\circ_i \Ga_2$ is defined by substituting the graph $\Ga_2$ into
the $i$-labeled vertex  $\xy
 (0,0)*+{_i}*\frm{o}\endxy$ of $\Ga_1$ and taking a sum over all possible re-attachments
of dangling half-edges (attached before to $v_i$) to the vertices of $\Ga_2$. For example, for $d$ odd one has
$$
\Ba{c}\resizebox{14mm}{!}{
\xy
%
   {\ar@{.>}@/^0.6pc/(-6,0)*+{_1}*\frm{o} ;(6,0)*+{_2}*\frm{o}};
 {\ar@{<.}@/^0.6pc/(6,0)*+{_2}*\frm{o};(-6,0)*+{_1}*\frm{o}};
\endxy}
\Ea
 \ \ \circ_1\
 \Ba{c} \resizebox{2.8mm}{!}{\xy
 (0,5)*+{_1}*\frm{o}="A";
 (0,-5)*+{_2}*\frm{o}="B";
 \ar @{.>} "A";"B" <0pt>
\endxy} \Ea
=
 \Ba{c}\resizebox{11mm}{!}{
\xy
 (-5,9)*+{_1}*\frm{o}="A";
 (-5,0)*+{_2}*\frm{o}="B";
 {\ar @{.>} "A";"B" <0pt>};
  {\ar@{.>}@/^0.6pc/(-5,0)*+{_2}*\frm{o} ;(5,0)*+{_3}*\frm{o}};
 {\ar@{<.}@/^0.6pc/(5,0)*+{_3}*\frm{o};(-5,0)*+{_2}*\frm{o}};
\endxy}\Ea
+
 \Ba{c}\resizebox{11mm}{!}{
\xy
 (-5,0)*+{_1}*\frm{o}="A";
 (-5,-9)*+{_2}*\frm{o}="B";
 {\ar @{.>} "A";"B" <0pt>};
  {\ar@{.>}@/^0.6pc/(-5,0)*+{_1}*\frm{o} ;(5,0)*+{_3}*\frm{o}};
 {\ar@{<.}@/^0.6pc/(5,0)*+{_3}*\frm{o};(-5,0)*+{_1}*\frm{o}};
\endxy}\Ea
+2
 \Ba{c}\resizebox{10mm}{!}{
 \xy
 (0,5)*+{_1}*\frm{o}="A";
 (0,-5)*+{_2}*\frm{o}="B";
 (10,0)*+{_3}*\frm{o}="C";
 {\ar @{.>} "A";"B" <0pt>};
 {\ar @{.>} "A";"C" <0pt>};
 {\ar @{.>} "B";"C" <0pt>};
\endxy}\Ea.
$$
Note that for $d$ even the first graph in the above formula vanishes identically as it has an automorphism reversing its orientation.

\subsubsection{\bf Operad of degree shifted Lie algebras}
The operad of degree  $d\in \Z$-shifted Lie algebras is the quotient
$$
\Lie_{d}:=\cF ree\langle E\rangle/I,
$$
of the free operad generated by an  $\bS$-module $E=\{E(n)\}_{n\geq 2}$
$$
E(n):=\left\{ \Ba{ll} \sgn_2^{\ot |d|}\ot \id_1[d-1]=\mbox{span}\left\langle
\Ba{c}\begin{xy}
 <0mm,0.66mm>*{};<0mm,3mm>*{}**@{-},
 <0.39mm,-0.39mm>*{};<2.2mm,-2.2mm>*{}**@{-},
 <-0.35mm,-0.35mm>*{};<-2.2mm,-2.2mm>*{}**@{-},
 <0mm,0mm>*{\bu};<0mm,0mm>*{}**@{},
   <0.39mm,-0.39mm>*{};<2.9mm,-4mm>*{^{_2}}**@{},
   <-0.35mm,-0.35mm>*{};<-2.8mm,-4mm>*{^{_1}}**@{},
\end{xy}\Ea
=(-1)^{d}
\Ba{c}\begin{xy}
 <0mm,0.66mm>*{};<0mm,3mm>*{}**@{-},
 <0.39mm,-0.39mm>*{};<2.2mm,-2.2mm>*{}**@{-},
 <-0.35mm,-0.35mm>*{};<-2.2mm,-2.2mm>*{}**@{-},
 <0mm,0mm>*{\bu};<0mm,0mm>*{}**@{},
   <0.39mm,-0.39mm>*{};<2.9mm,-4mm>*{^{_1}}**@{},
   <-0.35mm,-0.35mm>*{};<-2.8mm,-4mm>*{^{_2}}**@{},
\end{xy}\Ea
\right\rangle  & \text{if}\ n=2,\\
0 & \text{otherwise}.
\Ea\right.
$$
 modulo the ideal generated by the following relation\footnote{When representing elements of operads and props as decorated graphs we tacitly assume that all edges and legs are {\em directed}\, along the flow going from the bottom of the graph to the top}.
$$
\Ba{c}\resizebox{10mm}{!}{ \begin{xy}
 <0mm,0mm>*{\bu};<0mm,0mm>*{}**@{},
 <0mm,0.69mm>*{};<0mm,3.0mm>*{}**@{-},
 <0.39mm,-0.39mm>*{};<2.4mm,-2.4mm>*{}**@{-},
 <-0.35mm,-0.35mm>*{};<-1.9mm,-1.9mm>*{}**@{-},
 <-2.4mm,-2.4mm>*{\bu};<-2.4mm,-2.4mm>*{}**@{},
 <-2.0mm,-2.8mm>*{};<0mm,-4.9mm>*{}**@{-},
 <-2.8mm,-2.9mm>*{};<-4.7mm,-4.9mm>*{}**@{-},
    <0.39mm,-0.39mm>*{};<3.3mm,-4.0mm>*{^3}**@{},
    <-2.0mm,-2.8mm>*{};<0.5mm,-6.7mm>*{^2}**@{},
    <-2.8mm,-2.9mm>*{};<-5.2mm,-6.7mm>*{^1}**@{},
 \end{xy}}\Ea
 +
  \Ba{c}\resizebox{10mm}{!}{ \begin{xy}
 <0mm,0mm>*{\bu};<0mm,0mm>*{}**@{},
 <0mm,0.69mm>*{};<0mm,3.0mm>*{}**@{-},
 <0.39mm,-0.39mm>*{};<2.4mm,-2.4mm>*{}**@{-},
 <-0.35mm,-0.35mm>*{};<-1.9mm,-1.9mm>*{}**@{-},
 <-2.4mm,-2.4mm>*{\bu};<-2.4mm,-2.4mm>*{}**@{},
 <-2.0mm,-2.8mm>*{};<0mm,-4.9mm>*{}**@{-},
 <-2.8mm,-2.9mm>*{};<-4.7mm,-4.9mm>*{}**@{-},
    <0.39mm,-0.39mm>*{};<3.3mm,-4.0mm>*{^2}**@{},
    <-2.0mm,-2.8mm>*{};<0.5mm,-6.7mm>*{^1}**@{},
    <-2.8mm,-2.9mm>*{};<-5.2mm,-6.7mm>*{^3}**@{},
 \end{xy}}\Ea
 +
\Ba{c}\resizebox{10mm}{!}{ \begin{xy}
 <0mm,0mm>*{\bu};<0mm,0mm>*{}**@{},
 <0mm,0.69mm>*{};<0mm,3.0mm>*{}**@{-},
 <0.39mm,-0.39mm>*{};<2.4mm,-2.4mm>*{}**@{-},
 <-0.35mm,-0.35mm>*{};<-1.9mm,-1.9mm>*{}**@{-},
 <-2.4mm,-2.4mm>*{\bu};<-2.4mm,-2.4mm>*{}**@{},
 <-2.0mm,-2.8mm>*{};<0mm,-4.9mm>*{}**@{-},
 <-2.8mm,-2.9mm>*{};<-4.7mm,-4.9mm>*{}**@{-},
    <0.39mm,-0.39mm>*{};<3.3mm,-4.0mm>*{^3}**@{},
    <-2.0mm,-2.8mm>*{};<0.5mm,-6.7mm>*{^2}**@{},
    <-2.8mm,-2.9mm>*{};<-5.2mm,-6.7mm>*{^1}**@{},
 \end{xy}}\Ea
 =0
$$

\subsubsection{\bf Kontsevich graph complex as a deformation complex}
There is a morphism of operads \cite{Wi1}
$$
\Ba{rccc}
i: & \Lie_d & \lon &  \Gra_d\\
   & \Ba{c}\begin{xy}
 <0mm,0.66mm>*{};<0mm,3mm>*{}**@{-},
 <0.39mm,-0.39mm>*{};<2.2mm,-2.2mm>*{}**@{-},
 <-0.35mm,-0.35mm>*{};<-2.2mm,-2.2mm>*{}**@{-},
 <0mm,0mm>*{\bu};<0mm,0mm>*{}**@{},
   <0.39mm,-0.39mm>*{};<2.9mm,-4mm>*{^{_2}}**@{},
   <-0.35mm,-0.35mm>*{};<-2.8mm,-4mm>*{^{_1}}**@{},
\end{xy}\Ea
&\lon&
\Ba{c} \resizebox{11mm}{!}{\xy
 (10,1)*+{_2}*\frm{o}="B";
 (0,1)*+{_1}*\frm{o}="A";
 \ar @{.>} "A";"B" <0pt>
\endxy} \Ea
\Ea
$$
and the full Kontsevich graph complex can be defined as the dg Lie algebra,
$$
\fGC_d=\Def(\Lie_d \rar \Gra_d)\simeq \prod_{n\geq 1} \Gra(n)\ot_{\bS_n} \sgn_n^{\ot |d|}[d(1-n)]
$$
controlling deformations of this morphism. 
It is equipped with the pre-Lie algebra structure given by the composition
\Beq\label{2: preLie in fcGCd}
\Ga_1\circ \Ga_2:=\sum_{v\in V(\Ga_1)} \Ga_1\circ_v \Ga_2
\Eeq
where $\Ga_1\circ_v \Ga_2$ stands for the substitution of $\Ga_2$ into the vertex $v$ and redistribution of edges attached earlier to $v$ among vertices of $\Ga_2$ in all possible ways. Notice that the vertices of generators $\Ga$ of $\fGC_d$ are unlabelled, and in the case $d$
odd one also assumes that vertices are ordered (up to an odd permutation and multiplying
by $-1$). The cohomological degree of $\Ga$ is given by
$$
|\Ga|=d(\# V(\Ga) -1) + (1-d)\# E(\Ga)
$$
and the Lie bracket in $\fGC_d$ is given by
$$
[\Ga_1,\Ga_2]= \Ga_1\circ \Ga_2 - (-1)^{|\Ga_1||\Ga_2|} \Ga_2\circ \Ga_1.
$$
The morphism of operads $i$ gets encoded now into a Maurer-Cartan element of
$(\fGC_d,[\ ,\ ])$,
$$
[\ga,\ga]=0,
$$
given explicitly by
$$
\ga_0:=
\Ba{c}\resizebox{7mm}{!}{\xy
%
 (0,1)*{\circ}="a",
(8,1)*{\circ}="b",
\ar @{.>} "a";"b" <0pt>
\endxy}\Ea \in \fGC_d
$$
The differential in $\fGC_d$ is given by the commutator
with $\ga$,
$$
\delta \Ga:=[\ga_0,\Ga ].
$$

We consider the dg Lie subalgebra  
$$
\fcGC_d\subset \fGC_d
$$
spanned by connected graphs.
The dg Lie algebra $\fcGC_d$ contains another dg Lie subalgebra $\GC_d^{2}$ generated by graphs
with every vertex of valence $\geq 2$; moreover the inclusion
$$
\GC_d^{2}\lon \fcGC_d
$$
is a quasi-isomorphism of dg Lie algebras \cite{Wi1}.

\sip

The dg Lie algebra  $\GC_d^{2}$ contains a dg Lie subalgebra  $\GC_d$
generated by graphs with every vertex of valence $\geq 3$; the inclusion $\GC_d$ into
$\GC_d^{2}$ is not a quasi-isomorphism, but the difference at the cohomology level is well-understood \cite{Wi1},
$$
H^\bu(\GC_d^{2})= H^\bu(\GC_d)\ \oplus \
\bigoplus_{j\geq 1\atop j\equiv 2d+1 \mod 4} \K[d-j]
$$
where the summand $\K[d-j]$ is generated by the polytopes-like graphs with $j$ vertices. It is proven
in \cite{Wi1} that there is an isomorphism of Lie algebras,
$$
H^0(\GC_2)=\grt_1,
$$
where $\grt_1$ is the Lie algebra of the Grothendieck-Teichm\"uller group $GRT_1$.
Moreover, $H^{\bu<0}(\GC_2)=0$.

\subsection{Oriented graph complex} Let $\Gra_{d+1}^{or}$ be a version of $\Gra_{d+1}$ in which the edges of the generating graphs have a fixed direction, i.e.\ the rrelation (\ref{2: symmetry of dotted edges}) is dropped, but the directed edges are not allowed to form closed directed paths. We show edges of graphs from   $\Gra_d^{or}$ as {\em solid}\, (rather than dotted) arrows, while vertices are depicted in the black colour, i.e.
$$
\Ba{c}\resizebox{10mm}{!}{
\xy
(-5,2)*{^1},
(5,2)*{^2},
   {\ar@/^0.6pc/(-5,0)*{\bbu};(5,0)*{\bu}};
 {\ar@{<-}@/^0.6pc/(5,0)*{\bbu};(-5,0)*{\bu}};
\endxy}\Ea\in \Gra_{d+1}^{or} , \ \
\Ba{c}\resizebox{10mm}{!}{
\xy
(-5,2)*{^1},
(5,2)*{^2},
   {\ar@/^0.6pc/(-5,0)*{\bbu};(5,0)*{\bu}};
 {\ar@/^0.6pc/(5,0)*{\bbu};(-5,0)*{\bu}};
\endxy}
\Ea\not\in \Gra_{d+1}^{or}.
$$
There is a morphism of operads
$$
\Ba{rccc}
i: & \Lie_{d+1} & \lon &  \Gra_{d+1}^{or}\\
   & \Ba{c}\begin{xy}
 <0mm,0.66mm>*{};<0mm,3mm>*{}**@{-},
 <0.39mm,-0.39mm>*{};<2.2mm,-2.2mm>*{}**@{-},
 <-0.35mm,-0.35mm>*{};<-2.2mm,-2.2mm>*{}**@{-},
 <0mm,0mm>*{\bu};<0mm,0mm>*{}**@{},
   <0.39mm,-0.39mm>*{};<2.9mm,-4mm>*{^{_2}}**@{},
   <-0.35mm,-0.35mm>*{};<-2.8mm,-4mm>*{^{_1}}**@{},
\end{xy}\Ea
&\lon&
\frac{1}{2}\left(
\Ba{c}\resizebox{9.6mm}{!}{
\xy
 (0,0)*{\bbu}="a",
(9,0)*{\bbu}="b",
 (0,2)*{^1},
(9,2)*{^2},
\ar @{->} "a";"b" <0pt>
\endxy
} \Ea
-(-1)^{d+1}
\Ba{c}\resizebox{9.6mm}{!}{
\xy
 (0,0)*{\bbu}="a",
(9,0)*{\bbu}="b",
 (0,2)*{^1},
(9,2)*{^2},
\ar @{<-} "a";"b" <0pt>
\endxy
} \Ea\right)
\Ea
$$
and the {\em oriented}\, version of the full Kontsevich graph complex is defined as the dg Lie algebra,
$$
\fOGC_{d+1}=\Def(\Lie_{d+1} \rar \Gra_{d+1}^{or})
$$
controlling deformations of the above morphism. It contains a dg Lie subalgebra
$\OGC_{d+1}^0\subset \fOGC_{d+1}$ spanned by connected graphs which in turn contains a dg Lie subalgebra
$\OGC_{d+1}^{2}\subset \OGC_{d+1}^0$ spanned by graphs which have each vertex at least bivalent. The inclusion
$$
\OGC_{d+1}^{2} \lon \OGC_{d+1}^{0}
$$
is a quasi-isomorphism of dg Lie algebras. It has been proven in \cite{Wi2} that at the cohomology level there is an isomorphism of Lie algebras,
\Beq\label{2: GCd versus OGC_d+1}
H^\bu(\GC_d^{2})\simeq H^\bu(\OGC_{d+1}^{2})
\Eeq
for any $d\in \Z$. A different proof of the above isomorphism at the level of graded vector space was given by Marko \v Zivkovi\'c in \cite{Z1}
by constructing an {\em explicit}\, morphism on (slightly reduced versions of the) dual graph complexes, not respecting, however, the Lie co-algebra structures. One more proof of the isomorphism (\ref{2: GCd versus OGC_d+1})
at the level of graded vector spaces was given in \cite{Me0}. In this paper we improve further
 the first and the third approaches to show that the isomorphism  (\ref{2: GCd versus OGC_d+1}) lifts
to the chains level, i.e.\ that the Kontsevich graphs complex $\GC_d^{2}$ and its oriented version $\OGC_{d+1}^{2}$ are quasi-isomorphic as $\Lie_\infty$ algebras.

\bip

\bip


{\large
\section{\bf A 2-coloured operad of oriented graphs}
}
\label{3 sec:2-coloured graph operad}

\mip

\subsection{A 2-coloured operad of homotopy Lie algebras}
The minimal resolution of the operad $\Lie_{d+1}$ (see \S {\bf 2.1.1})
 is a dg free operad whose (skew)symmetric generators,
\Beq\label{2: Lie_inf corolla}
\Ba{c}\resizebox{22mm}{!}{ \xy
(1,-5)*{\ldots},
(-13,-7)*{_1},
(-8,-7)*{_2},
(-3,-7)*{_3},
(7,-7)*{_{n-1}},
(13,-7)*{_n},
 (0,0)*{\bu}="a",
(0,5)*{}="0",
(-12,-5)*{}="b_1",
(-8,-5)*{}="b_2",
(-3,-5)*{}="b_3",
(8,-5)*{}="b_4",
(12,-5)*{}="b_5",
\ar @{-} "a";"0" <0pt>
\ar @{-} "a";"b_2" <0pt>
\ar @{-} "a";"b_3" <0pt>
\ar @{-} "a";"b_1" <0pt>
\ar @{-} "a";"b_4" <0pt>
\ar @{-} "a";"b_5" <0pt>
\endxy}\Ea
=(-1)^{d+1}
\Ba{c}\resizebox{23mm}{!}{\xy
(1,-6)*{\ldots},
(-13,-7)*{_{\sigma(1)}},
(-6.7,-7)*{_{\sigma(2)}},
(13,-7)*{_{\sigma(n)}},
 (0,0)*{\bu}="a",
(0,5)*{}="0",
(-12,-5)*{}="b_1",
(-8,-5)*{}="b_2",
(-3,-5)*{}="b_3",
(8,-5)*{}="b_4",
(12,-5)*{}="b_5",
\ar @{-} "a";"0" <0pt>
\ar @{-} "a";"b_2" <0pt>
\ar @{-} "a";"b_3" <0pt>
\ar @{-} "a";"b_1" <0pt>
\ar @{-} "a";"b_4" <0pt>
\ar @{-} "a";"b_5" <0pt>
\endxy}\Ea,
\ \ \ \forall \sigma\in \bS_n,\ n\geq2,
\Eeq
have degrees $2+d-n(d+1)$; the case $d=0$ corresponds to the usual operad of Lie algebras
(when Lie brackets have degree zero) and to its minimal resolution.
 The differential in $\Holie_{d+1}$ is given by
\Beq\label{3: d in Lie_infty}
\delta\hspace{-3mm}
\Ba{c}\resizebox{21mm}{!}{\xy
(1,-5)*{\ldots},
(-13,-7)*{_1},
(-8,-7)*{_2},
(-3,-7)*{_3},
(7,-7)*{_{n-1}},
(13,-7)*{_n},
 (0,0)*{\bu}="a",
(0,5)*{}="0",
(-12,-5)*{}="b_1",
(-8,-5)*{}="b_2",
(-3,-5)*{}="b_3",
(8,-5)*{}="b_4",
(12,-5)*{}="b_5",
\ar @{-} "a";"0" <0pt>
\ar @{-} "a";"b_2" <0pt>
\ar @{-} "a";"b_3" <0pt>
\ar @{-} "a";"b_1" <0pt>
\ar @{-} "a";"b_4" <0pt>
\ar @{-} "a";"b_5" <0pt>
\endxy}\Ea
=
\sum_{A\varsubsetneq [n]\atop
\# A\geq 2}\pm
%
%
\Ba{c}\resizebox{19mm}{!}{\begin{xy}
<10mm,0mm>*{\bu},
<10mm,0.8mm>*{};<10mm,5mm>*{}**@{-},
<0mm,-10mm>*{...},
<14mm,-5mm>*{\ldots},
<13mm,-7mm>*{\underbrace{\ \ \ \ \ \ \ \ \ \ \ \ \  }},
<14mm,-10mm>*{_{[n]\setminus A}};
<10.3mm,0.1mm>*{};<20mm,-5mm>*{}**@{-},
<9.7mm,-0.5mm>*{};<6mm,-5mm>*{}**@{-},
<9.9mm,-0.5mm>*{};<10mm,-5mm>*{}**@{-},
<9.6mm,0.1mm>*{};<0mm,-4.4mm>*{}**@{-},
<0mm,-5mm>*{\bu};
<-5mm,-10mm>*{}**@{-},
<-2.7mm,-10mm>*{}**@{-},
<2.7mm,-10mm>*{}**@{-},
<5mm,-10mm>*{}**@{-},
<0mm,-12mm>*{\underbrace{\ \ \ \ \ \ \ \ \ \ }},
<0mm,-15mm>*{_{A}}.
\end{xy}}
\Ea
\Eeq
If $d$ is odd, all the signs above are equal to $-1$.

\sip

If we relax the restriction $n\geq 2$ on the generators above to  $n \geq 1$, and the restriction
in the definition of the differential from  $\#A\geq 2$ to $\#A \geq 1$, then one gets a dg free operad
$\Holie_{d+1}^+$ which is acyclic but is often very helpful in building important deformation complexes.
The newly added generator
$\Ba{c}\resizebox{1.8mm}{!}{\begin{xy}
 <0mm,-0.55mm>*{};<0mm,-3mm>*{}**@{-},
 <0mm,0.5mm>*{};<0mm,3mm>*{}**@{-},
 <0mm,0mm>*{\bu};<0mm,0mm>*{}**@{},
 \end{xy}}\Ea
$  controls deformations of the differential in representations
of $\Holie_{d+1}^+$ in dg vector spaces.

\sip

Let $\Holie_{d,d+1}^{+}$ be the 2-coloured operad controlling $\Holie_1^+$-algebra structures in the
 dg vector spaces which are direct sums
$$
V[d]\oplus W[d-1]
$$
for arbitrary dg spaces $V$ and $W$; put another way, a representation of $\Holie_{d,d+1}^{+}$ in the
2-coloured endomorphism operad $\cE nd_{V,W}$ is the same as a
$\Holie_1^+$-algebra structure in the above direct sum. Thus  $\Holie_{d,d+1}^{+}$ is a
dg free 2-coloured operad generated by two sets of corollas, one set having the output in, say, white colour
(shown in pictures as the dotted legs), and the other set having an output in black colour 
(shown as solid legs)
$$
\fC_{m+n}^\circ:= \Ba{c}\resizebox{31mm}{!}{
\begin{xy}
 <0mm,-0.5mm>*{\bu};
 <0mm,0mm>*{};<0mm,5mm>*{}**@{.},
 <0mm,0mm>*{};<-16mm,-5mm>*{}**@{-},
 <0mm,0mm>*{};<-11mm,-5mm>*{}**@{-},
 <0mm,0mm>*{};<-3.5mm,-5mm>*{}**@{-},
 <0mm,0mm>*{};<-6mm,-5mm>*{...}**@{},
   <0mm,0mm>*{};<-16mm,-8mm>*{^{1}}**@{},
   <0mm,0mm>*{};<-11mm,-8mm>*{^{2}}**@{},
   <0mm,0mm>*{};<-3mm,-8mm>*{^{m}}**@{},
 <0mm,0mm>*{};<16mm,-5mm>*{}**@{.},
 <0mm,0mm>*{};<8mm,-5mm>*{}**@{.},
 <0mm,0mm>*{};<3.5mm,-5mm>*{}**@{.},
 <0mm,0mm>*{};<11.6mm,-5mm>*{...}**@{},
   <0mm,0mm>*{};<19mm,-8mm>*{^{\bar{n}}}**@{},
<0mm,0mm>*{};<10mm,-8mm>*{^{\bar{2}}}**@{},
   <0mm,0mm>*{};<5mm,-8mm>*{^{\bar{1}}}**@{},
 \end{xy}}
\Ea\in sgn_m^{\ot |d+1|}\ot sgn_n^{|d|}[(d+1)m+ dn -2-d]
, \ \ \ \ m+n\geq 1,
$$
$$
\fC_{m+n}^\bu=  \Ba{c}\resizebox{31mm}{!}{
\begin{xy}
 <0mm,-0.5mm>*{\bu};
 <0mm,0mm>*{};<0mm,5mm>*{}**@{-},
 <0mm,0mm>*{};<-16mm,-5mm>*{}**@{-},
 <0mm,0mm>*{};<-11mm,-5mm>*{}**@{-},
 <0mm,0mm>*{};<-3.5mm,-5mm>*{}**@{-},
 <0mm,0mm>*{};<-6mm,-5mm>*{...}**@{},
   <0mm,0mm>*{};<-18mm,-8mm>*{^{1}}**@{},
   <0mm,0mm>*{};<-11mm,-8mm>*{^{2}}**@{},
   <0mm,0mm>*{};<-3mm,-8mm>*{^{m}}**@{},
 <0mm,0mm>*{};<16mm,-5mm>*{}**@{.},
 <0mm,0mm>*{};<8mm,-5mm>*{}**@{.},
 <0mm,0mm>*{};<3.5mm,-5mm>*{}**@{.},
 <0mm,0mm>*{};<11.6mm,-5mm>*{...}**@{},
   <0mm,0mm>*{};<19mm,-8mm>*{^{\bar{n}}}**@{},
<0mm,0mm>*{};<10mm,-8mm>*{^{\bar{2}}}**@{},
   <0mm,0mm>*{};<5mm,-8mm>*{^{\bar{1}}}**@{},
 \end{xy}}
\Ea \in sgn_m^{\ot |d+1|}\ot sgn_n^{|d|}[(d+1)m+ dn -1-d]
, \ \ \ \ m+n\geq 1.
$$
We understand these corollas as generators of the 1-dimensional $\bS_m\times \bS_n$ modules
shown above. The differential in $\Holie_{d,d+1}^{+}$ is given by splitting each corolla
into two sums, one is via the substitution of the dotted edge inside the unique vertex,
and the other is via the substituting the solid edge (cf.\ (\ref{3: d in Lie_infty})).

\sip

Let $I$ be the ideal in $\Holie_{d,d+1}^{+}$ generated by corollas $\fC^\circ_{1+0}$,
$\fC_{1+0}$, and all corollas  $\fC_{m+n}^\bu$
with $m\geq 1$. It is closed under the differential so that
the quotient 2-coloured operad
$$
\Holie_{d,d+1}:= \Holie_{d,d+1}^+/I
$$
is a dg {\em free}\, 2-colored operad with the following generators,
$$
\left\{
\Ba{c}\resizebox{31mm}{!}{
\begin{xy}
 <0mm,-0.5mm>*{\bu};
 <0mm,0mm>*{};<0mm,5mm>*{}**@{.},
 <0mm,0mm>*{};<-16mm,-5mm>*{}**@{-},
 <0mm,0mm>*{};<-11mm,-5mm>*{}**@{-},
 <0mm,0mm>*{};<-3.5mm,-5mm>*{}**@{-},
 <0mm,0mm>*{};<-6mm,-5mm>*{...}**@{},
   <0mm,0mm>*{};<-16mm,-8mm>*{^{1}}**@{},
   <0mm,0mm>*{};<-11mm,-8mm>*{^{2}}**@{},
   <0mm,0mm>*{};<-3mm,-8mm>*{^{m}}**@{},
 <0mm,0mm>*{};<16mm,-5mm>*{}**@{.},
 <0mm,0mm>*{};<8mm,-5mm>*{}**@{.},
 <0mm,0mm>*{};<3.5mm,-5mm>*{}**@{.},
 <0mm,0mm>*{};<11.6mm,-5mm>*{...}**@{},
   <0mm,0mm>*{};<19mm,-8mm>*{^{\bar{n}}}**@{},
<0mm,0mm>*{};<10mm,-8mm>*{^{\bar{2}}}**@{},
   <0mm,0mm>*{};<5mm,-8mm>*{^{\bar{1}}}**@{},
 \end{xy}}
\Ea m+n\geq 2,
\ \ \
 \Ba{c}\resizebox{1.8mm}{!}{\begin{xy}
 <0mm,-0.55mm>*{};<0mm,-3mm>*{}**@{-},
 <0mm,0.5mm>*{};<0mm,3mm>*{}**@{.},
 <0mm,0mm>*{\bu};<0mm,0mm>*{}**@{},
 \end{xy}}\Ea,
 \ \ \
 \Ba{c}\resizebox{22mm}{!}{ \xy
(1,-5)*{\ldots},
(-13,-7)*{_1},
(-8,-7)*{_2},
(-3,-7)*{_3},
(7,-7)*{_{n-1}},
(13,-7)*{_n},
 (0,0)*{\bu}="a",
(0,5)*{}="0",
(-12,-5)*{}="b_1",
(-8,-5)*{}="b_2",
(-3,-5)*{}="b_3",
(8,-5)*{}="b_4",
(12,-5)*{}="b_5",
\ar @{-} "a";"0" <0pt>
\ar @{-} "a";"b_2" <0pt>
\ar @{-} "a";"b_3" <0pt>
\ar @{-} "a";"b_1" <0pt>
\ar @{-} "a";"b_4" <0pt>
\ar @{-} "a";"b_5" <0pt>
\endxy}\Ea n\geq 2
\right\}
$$
and with the differential given by (\ref{3: d in Lie_infty}) on generators with black output, and by the following formula
$$
\delta
\Ba{c}\resizebox{32mm}{!}{
\begin{xy}
 <0mm,-0.5mm>*{\bu};
 <0mm,0mm>*{};<0mm,5mm>*{}**@{.},
 <0mm,0mm>*{};<-16mm,-5mm>*{}**@{-},
 <0mm,0mm>*{};<-11mm,-5mm>*{}**@{-},
 <0mm,0mm>*{};<-3.5mm,-5mm>*{}**@{-},
 <0mm,0mm>*{};<-6mm,-5mm>*{...}**@{},
   <0mm,0mm>*{};<-16mm,-8mm>*{^{1}}**@{},
   <0mm,0mm>*{};<-11mm,-8mm>*{^{2}}**@{},
   <0mm,0mm>*{};<-3mm,-8mm>*{^{m}}**@{},
 <0mm,0mm>*{};<16mm,-5mm>*{}**@{.},
 <0mm,0mm>*{};<8mm,-5mm>*{}**@{.},
 <0mm,0mm>*{};<3.5mm,-5mm>*{}**@{.},
 <0mm,0mm>*{};<11.6mm,-5mm>*{...}**@{},
   <0mm,0mm>*{};<19mm,-8mm>*{^{\bar{n}}}**@{},
<0mm,0mm>*{};<10mm,-8mm>*{^{\bar{2}}}**@{},
   <0mm,0mm>*{};<5mm,-8mm>*{^{\bar{1}}}**@{},
 \end{xy}}
\Ea
=
\sum_{A\varsubsetneq [n]\atop
\# A\geq 2}\ \pm\hspace{-3mm}
\Ba{c}\resizebox{39mm}{!}{
\begin{xy}
 <0mm,-0.5mm>*{\bu};
 <0mm,0mm>*{};<0mm,5mm>*{}**@{.},
 <0mm,0mm>*{};<-16mm,-5mm>*{}**@{-},
 <0mm,0mm>*{};<-11mm,-5mm>*{}**@{-},
 <0mm,0mm>*{};<-3.5mm,-5mm>*{}**@{-},
 <0mm,0mm>*{};<-6mm,-5mm>*{...}**@{},
 <0mm,0mm>*{};<16mm,-5mm>*{}**@{.},
 <0mm,0mm>*{};<8mm,-5mm>*{}**@{.},
 <0mm,0mm>*{};<3.5mm,-5mm>*{}**@{.},
 <0mm,0mm>*{};<11.6mm,-5mm>*{...}**@{},
   <0mm,0mm>*{};<17mm,-8mm>*{^{\bar{m}}}**@{},
<0mm,0mm>*{};<10mm,-8mm>*{^{\bar{2}}}**@{},
   <0mm,0mm>*{};<5mm,-8mm>*{^{\bar{1}}}**@{},
<-17mm,-12mm>*{\underbrace{\ \ \ \ \ \ \ \ \ \   }},
<-17mm,-14.9mm>*{_A};
<-6mm,-7mm>*{\underbrace{\ \ \ \ \ \ \  }},
<-6mm,-10mm>*{_{[n]\setminus A}};
 (-16.5,-5.5)*{\bu}="a",
(-23,-10)*{}="b_1",
(-20,-10)*{}="b_2",
(-16,-10)*{...}="b_3",
(-12,-10)*{}="b_4",
\ar @{-} "a";"b_2" <0pt>
\ar @{-} "a";"b_1" <0pt>
\ar @{-} "a";"b_4" <0pt>
 \end{xy}}
\Ea
+
 \sum_{ [m]=I_1\sqcup I_2\atop
 [n]=J_1\sqcup J_2}\hspace{-2mm}\pm
\Ba{c}
\begin{xy}
 <0mm,-0.5mm>*{\bu};
 <0mm,0mm>*{};<0mm,6mm>*{}**@{.},
 <0mm,0mm>*{};<-16mm,-6mm>*{}**@{-},
 <0mm,0mm>*{};<-11mm,-6mm>*{}**@{-},
 <0mm,0mm>*{};<-3.5mm,-6mm>*{}**@{-},
 <0mm,0mm>*{};<-6mm,-6mm>*{...}**@{},
%
%
 <0mm,0mm>*{};<3mm,-7mm>*{}**@{.},
 <0mm,0mm>*{};<8mm,-6mm>*{}**@{.},
 <0mm,0mm>*{};<13mm,-6mm>*{}**@{.},
<0mm,0mm>*{};<10mm,-6mm>*{...}**@{},
<-1.5mm,-16mm>*{\underbrace{\ \ \ \  }_{I_2}},
<9mm,-16mm>*{\underbrace{\ \ \   }_{J_1}},
<11mm,-9mm>*{\underbrace{\    }_{J_2}},
<-10mm,-9mm>*{\underbrace{\ \ \ \ \ \ \ \ \ \ \ \   }_{I_1}},
 %
 (3,-7)*{\bu}="a",
(-4,-13)*{}="b_1",
(1,-13)*{}="b_2",
(8,-13)*{...},
(-1,-13)*{...},
(5,-13)*{}="b_3",
(11,-13)*{}="b_4",
\ar @{-} "a";"b_2" <0pt>
\ar @{.} "a";"b_3" <0pt>
\ar @{-} "a";"b_1" <0pt>
\ar @{.} "a";"b_4" <0pt>
 \end{xy}
\Ea
$$
one the generators with white output. 
We  use the operad $\Holie_{d,d+1}$ together 
with a certain 2-coloured operad of graphs to build a new graph complex $\hOGC_{d,d+1}$
 below,  the main gadget which solves the long standing problem as stated
in the Main Theorem.

\subsection{A 2-coloured dg operad of graphs $\Gra_{d,d+1}$}

Let $G_{m,n;p}$ be the set of connected  graphs $\Ga$ with $m$
labelled (by integers from $\{1,\ldots,m\}$) white vertices,   $n$ labelled (by integers
from $\{\bar{1},\ldots, \bar{n}\}$)  black vertices, and $p$ labelled directed solid edges,  e.g.
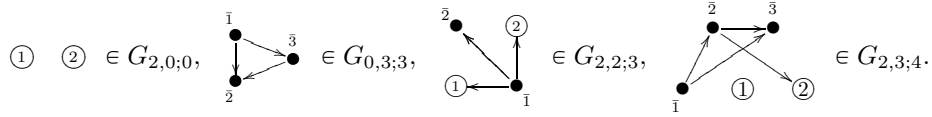
 $$
   \Ba{c}\resizebox{11mm}{!}{ \xy
(-4,1)*+{_1}*\frm{o}="1";
 (4,1)*+{_2}*\frm{o}="2";
\endxy} \Ea \in G_{2,0; 0}
,
 \Ba{c}\resizebox{11mm}{!}{
\xy
(-6,6.5)*{^{\bar{1}}},
(-6,-6.5)*{_{\bar{2}}},
(5,2.7)*{^{\bar{3}}},
 {\ar@{->}(-5,4)*{\bbu};(-5,-4)*{\bbu}};
  {\ar@{->}(-5,4)*{\bbu};(5,0)*{\bbu}};
 {\ar@{->}(5,0)*{\bbu};(-5,-4)*{\bbu}};
\endxy}\Ea\in G_{0,3;3}
,
\Ba{c}\resizebox{13mm}{!}{ \xy
(6.7,-2.3)*{_{\bar{1}}},
(-6.7,11.3)*{^{\bar{2}}},
(-5,0)*+{_1}*\frm{o}="1";
 (5,0)*{\bbu}="11";
 (5,10)*+{_2}*\frm{o}="2";
 (-5,10)*{\bbu}="22";
 \ar @{->} "11";"22" <0pt>
 \ar @{->} "11";"22" <0pt>
  \ar @{->} "11";"1" <0pt>
   \ar @{->} "11";"2" <0pt>
\endxy} \Ea \in G_{2,2; 3}
,
\Ba{c}\resizebox{20mm}{!}{ \xy
(-11,-2.8)*{_{\bar{1}}},
(-5,12.7)*{^{\bar{2}}},
(5,12.7)*{^{\bar{3}}},
(-10,0)*{\bbu}="11";
(0,0)*+{1}*\frm{o}="1";
(10,0)*+{2}*\frm{o}="2";
  (5,10)*{\bbu}="22";
 (-5,10)*{\bbu}="33";
 \ar @{->} "11";"33" <0pt>
 \ar @{->} "33";"22" <0pt>
 \ar @{<-} "22";"11" <0pt>
 \ar @{->} "33";"2" <0pt>
\endxy} \Ea \in G_{2,3;4}.
$$
We also assume that the generators $\Ga$
\Bi
\item[(a)] have no closed paths of directed edges, and
\item[(b)] have no {\em outgoing} edges attached to white vertices.
\Ei
 Let us call a {\em sink}\, any vertex which has no outgoing edges 
 attached; then every white vertex of a graph from $G_{m\geq 1,n;p,q}$ is a sink 
 (we put no such restriction on black vertices).

\sip

If we denote by $E(\Ga)$ the set of edges of $\Ga\in G_{m,n;p}$, then the {\em labelling}\, of
edges means that some isomorphism
$$
E(\Ga)\rar [p],
$$
is fixed; we do not show these labellings in our pictures as soon we 
shall get rid of this
 extra datum.

\sip

The group $\bS_p$ acts on the $G_{m,n;p}$ by permuting labels of edges. Hence it makes sense to consider,  for any integer $d\in \Z$
and any fixed pair of natural numbers $m,n\in \N$, a dg vector space 
(in fact, a dg $\bS_m\times \bS_n$-module),
$$
\Gra_{d,d+1}^\circ(m,n):= \prod_{p\geq 0} \K \langle G_{m,n;p}\rangle
\ot_{ \bS_p} \sgn_p^{|d|}[pd]
$$
equipped with the zero differential.

One can understand a generator of the $\bS_m\times\bS_n$-module 
as a graph $\Ga$ with
\Bi
\item[(i)]  $m$ labeled
white vertices,
\item[(ii)]
$n$ labelled black vertices,
\item[(iii)]  some number of unlabeled edges whose directions are {\em fixed},
\item[(iv)] a choice (up to sign) of a unital basis vector (called the
{\em orientation}) of the following 1-dimensional Euclidean vector space
 $$
   or(\Ga)\in \left\{\Ba{ll} \det(E(\Ga))   & \text{if $d$ is even} \\
         \K  &    \text{if $d$ is odd}
         \Ea
   \right.
 $$
 that is an ordering of edges for $d$ even.
 \Ei
A generator $\Ga\in  \Gra_{d+1}(m,n)$ is assigned the cohomological degree
$$
  |\Ga|=  - d \# E_{sol}(\Ga),
$$
i.e.\ every edge is assigned the degree $-d$.

\sip

Next we define a dg 2-coloured (with colours called {\em white}\, and {\em black}) 
operad of graphs
$$
\Gra_{d,d+1}=\left\{ \Gra_{d,d+1}(p):=\Gra_{d,d+1}^\circ(p) \oplus 
\Gra_{d,d+1}^\bu(p)\right\}_{p\geq 1}
$$
where $\Gra_{d,d+1}^\circ(p)$ (resp., $\Gra_{d,d+1}^\bu(p)$) stands for the dg $\bS_p$-module spanned by elements with so called white (resp.\ black) output; 
more precisely,
$$
\Gra_{d,d+1}^\circ(p):=\bigoplus_{p=m+n} 
\text{Ind}^{\bS_p}_{\bS_m\times \bS_n} \Gra^\circ_{d,d+1}(m,n), \ \ \ \Gra_{d,d+1}^\bu(p):=\Gra_{d+1}^{or}(p).
$$
Thus a generator of $\Gra_{d,d+1}$ is a pair of graphs
$$
\left(\sGa^\circ\in  \Gra_{d,d+1}(m,n)^\circ, \ \sGa^\bu \in \Gra_{d+1}^{or}   \right)
$$
and the operadic composition,  $\forall\ i\in [p]$,
$$
\Ba{rccc}
\circ_i: &  \cG ra_{d,d+1}(p) \times \cG ra_{d,d+1}(q) &\lon & \cG ra_{d,d+1}(p+q-1),
 \\
         &       \left(\sGa_1=(\sGa_1^\circ,\sGa_1^\bu), 
         \sGa_2=(\sGa_2^\circ, \sGa_2^\bu)\right) &\lon &      \sGa_1\circ_i \sGa_2,
\Ea
$$
is given by
$$
\sGa_1\circ_i \sGa_2:=\left\{\Ba{ll} (\sGa_1^\circ \circ_i\sGa_2^\circ, 0)   &
\text{if $i$ is a white vertex in $\Ga_1^\circ$}\\
(\sGa_1^\circ \circ_i \sGa_2^\bu, 0)   & \text{if $i$ is a black vertex 
in $\sGa_1^\circ$}\\
(0, \sGa_1^\bu \circ_i \sGa_2^\bu)   & \text{if $i$ is a black vertex in $\sGa_1^\bu$}\\
\Ea
\right.
$$
The symbol $A \circ_i B$ in the right hand side of the above equality stands for
the substitution of a graph $B$ into the $i$-labelled vertex $v$ of a graph $A$ and  
taking a sum over all possible re-attachments of dangling edges 
(attached before to $v$) to the vertices of $B$ (cf.\ \S {\ref{2: subsec on Gra_d}}).

\subsection{Proposition}{\em There is a morphism of dg 2-coloured operads
$$
f: \Holie_{d,d+1} \lon \Gra_{d,d+1}
$$
which vanishes on all generators of $\Holie_{d,d+1}$ except the following ones,

}
$$
f:\left\{\Ba{ccc}
 \Ba{c}\begin{xy}
 <0mm,0.66mm>*{};<0mm,3mm>*{}**@{-},
 <0.39mm,-0.39mm>*{};<2.2mm,-2.2mm>*{}**@{-},
 <-0.35mm,-0.35mm>*{};<-2.2mm,-2.2mm>*{}**@{-},
 <0mm,0mm>*{\bu};<0mm,0mm>*{}**@{},
   <0.39mm,-0.39mm>*{};<2.9mm,-4mm>*{^{_2}}**@{},
   <-0.35mm,-0.35mm>*{};<-2.8mm,-4mm>*{^{_1}}**@{},
\end{xy}\Ea
&\lon&
\frac{1}{2}\left(
\Ba{c}\resizebox{9.6mm}{!}{
\xy
 (0,0)*{\bbu}="a",
(9,0)*{\bbu}="b",
 (0,2)*{^1},
(9,2)*{^2},
\ar @{->} "a";"b" <0pt>
\endxy
} \Ea
-(-1)^{d+1}
\Ba{c}\resizebox{9.6mm}{!}{
\xy
 (0,0)*{\bbu}="a",
(9,0)*{\bbu}="b",
 (0,2)*{^1},
(9,2)*{^2},
\ar @{<-} "a";"b" <0pt>
\endxy
} \Ea\right)\\
\Ba{c}\resizebox{24mm}{!}{
\begin{xy}
 <0mm,-0.5mm>*{\bu};
 <0mm,0mm>*{};<0mm,5mm>*{}**@{.},
 <0mm,0mm>*{};<-3.5mm,-5mm>*{}**@{-},
   <0mm,0mm>*{};<-3.9mm,-8mm>*{^{\bar{1}}}**@{},
 <0mm,0mm>*{};<16mm,-5mm>*{}**@{.},
 <0mm,0mm>*{};<8mm,-5mm>*{}**@{.},
 <0mm,0mm>*{};<3.5mm,-5mm>*{}**@{.},
 <0mm,0mm>*{};<11.6mm,-5mm>*{...}**@{},
   <0mm,0mm>*{};<19mm,-8mm>*{^{{n}}}**@{},
<0mm,0mm>*{};<10mm,-8mm>*{^{{2}}}**@{},
   <0mm,0mm>*{};<5mm,-8mm>*{^{{1}}}**@{},
 \end{xy}}
\Ea &\lon &
\Ba{c}\resizebox{19mm}{!}{
\xy
 (-8,-1)*+{_1}*\frm{o}="a1";
 (-3,-1)*+{_2}*\frm{o}="a2"; 
 (2,-1)*{...};
 (8,-1)*+{_n}*\frm{o}="a3";  
   (0,8)*{\bbu}="b",
   (0,10)*{^{\bar{1}}},
\ar @{<-} "a1";"b" <0pt>
\ar @{<-} "a2";"b" <0pt>
\ar @{<-} "a3";"b" <0pt>
\endxy}\Ea \\
\Ea
\right.
$$

This Proposition can be proven in two ways. The first poof goes via a 
direct checking that the above formulae respect the differentials in both sides
of the morphism $f$ (which is a straightforward but a bit tedious calculation).

\sip

We choose another approach which is more suitable for our purposes in the next section, 
an approach which starts with a consideration of the dg Lie algebra 
\Beq\label{3: def complex of zero morphism}
\fhOGC_{d,d+1}:=\Def\left(\Holie_{d,d+1}\stackrel{0}{\lon} \Gra_{d,d+1}\right)_{conn}
\Eeq
controlling deformations of the {\em zero}\, morphism of the aforementioned operads; 
the differential in $\fhOGC_{d,d+1}$ is zero (for now).
We in particular consider the connected graded Lie subalgebra
\[
\hOGC_{d,d+1}\subset \fhOGC_{d,d+1}.
\]
In view of the operadic composition formulae in $\Gra_{d,d+1}$ described above,
we conclude that the dg Lie algebra $\hOGC_{d,d+1}$ decomposes into a semidirect product,
\begin{equation}\label{equ:fcGCsplit}
\hOGC_{d,d+1}=\fcGCdd \rtimes \fcOGC_{d+1},
\end{equation}
 of the dg  Lie algebra  $\fcOGC_{d+1}$ studied above and the dg Lie algebra $\fcGCdd$
generated by graphs $\Ga^\circ$ with two types of unlabelled vertices,
black and white ones.
 $$
\Ba{c}\resizebox{14mm}{!}{
\xy
 (0,0)*{\bullet}="a",
(0,8)*{\bullet}="b",
(-7.5,-4.5)*{\circ}="c",
(7.5,-4.5)*{\circ}="d",
\ar @{<-} "a";"b" <0pt>
\ar @{->} "a";"c" <0pt>
\ar @{->} "b";"c" <0pt>
\ar @{->} "b";"d" <0pt>
\ar @{<-} "d";"a" <0pt>
\endxy}
\Ea
 \in \fcGCdd.
$$
The cohomological degree is given by the formula 
$$
|\Ga^\circ|=d\#V_\circ(\Ga) + (d+1)\# V_\bu(\Ga^\circ)
-d \# E(\Ga^\circ)  -d.
$$
The graded Lie algebra $\hOGC_{d,d+1}$ has an underlying dg pre-Lie 
algebra structure  which is fully analogous to the one
(\ref{2: preLie in fcGCd}) in $\fcGC_d$.
If we represent a generic element $\Ga$ of the graded Lie algebra $\hOGC_{d,d+1}$ as a pair
$\Ga=(\Ga^\circ,\Ga^\bu)$, with $\Ga^\circ\in \OGC_{d,d+1}$ and $\Ga^\bu\in \OGC_{d+1}^0$, then the pre-Lie composition is given by
$$
(\Ga_1^\circ,\Ga_1^\bu)\circ (\Ga_2^\circ,\Ga_2^\bu)
= \left(\sum_{v\in V_\circ} \Ga_1^\circ \circ_v  \Ga_2^\circ +  
\sum_{v\in V_\bu}\Ga_1^\circ \circ_v  \Ga_2^\bu, \
\sum_{v\in V(\Ga_1^\bu)}  \Ga_1^\bu \circ_v  \Ga_2^\bu\right)
$$
The calculations given in \S 5.1.1 and \S 5.2.2 of \cite{Me0} imply that the degree 1 element
$$
\ga=\left( \ga^\circ:=
\ \bu
+ \sum_{k=1}^\infty \frac{1}{k!}
\resizebox{14mm}{!}{
\xy
 (-5,-1)*{\circ}="a1",
  (5,-1)*{\circ}="a2",
   (-2,-1)*{\circ}="a3",
    (2,-1)*{...},
   (0,7)*{\bu}="b",
(0,-4)*{\underbrace{\hspace{12mm}}_k},
\ar @{<-} "a1";"b" <0pt>
\ar @{<-} "a2";"b" <0pt>
\ar @{<-} "a3";"b" <0pt>
\endxy}, \ \ \ \ \ga^\bu:=\resizebox{3mm}{!}{
\xy
 (0,6)*{\bullet}="a",
(0,-2)*{\bu}="b",
\ar @{->} "a";"b" <0pt>
\endxy}\ \
\right)
$$
in $\hOGC_{d,d+1}$ satisfies the equation,
$$
\ga\circ \ga=0,
$$
i.e.\ that $\Ga$ is a Maurer-Cartan element of the deformation complex 
(\ref{3: def complex of zero morphism}). Any such a MC element corresponds to some
morphism from $\Holie_{d,d+1}$ to $\Gra_{d,d+1}$, and the one shown just above corresponds
precisely to the morphism $f$. The proposition is proven.

\sip

The twisted dg Lie algebra  $(\hOGC_{d,d+1}, [\ ,\ ], \delta:= [\ga,\ ])$ 
controls the deformation theory of the morphism  $f$.
From now one we abbreviate this structure to $\hOGC_{d,d+1}$; it provides us with the main tool which we use  to prove the Main Theorem in the next section.


\mip

{\large
\section{\bf Proof of the Main Theorem}
}
We claim that there are morphisms of dg Lie algebras 
$$
\fcGC_d    \stackrel{\pi_1}{\longleftarrow}   \hOGC_{d,d+1}  \stackrel{\pi_2}{\lon}
\OGC_{d+1}^{0}.
$$
Here the right-hand arrow $\pi_2$ is the obvious projection onto the second factor in \eqref{equ:fcGCsplit}.

To construct the left-hand arrow $\pi_1$ we introduce the following notation.
Let us call a vertex of a graph $\Gamma\in \fcGCdd$ (or in $\OGC_{d+1}^0$) \emph{inessential} if it is of valence 2 with two outputs, and \emph{essential} otherwise.
Then we define $\pi_1$ as follows:
\begin{itemize}
    \item $\pi_1(\OGC_{d+1}^{0}) = 0$.
    \item For $\Gamma\in \fcGCdd$ with at least one essential black vertex we set $\pi_1(\Gamma)=0$.
    \item Suppose $\Gamma\in \fcGCdd$ has only inessential black vertices. Then necessarily each such vertex has two white neighbors.
    We may hence build a graph $\pi_1(\Gamma)\in \fcGC_d$ by retaining only the white vertices, and adding one edge between white vertices $(u,v)$ for every inessential vertex connected to $u$ and $v$. (Note that here possibly $u=v$.)
\[
\Ba{c}\resizebox{18mm}{!}{
\xy
%
 (0,0)*{\circ}="a",
(8,0)*{\bullet}="b",
(16,0)*{\circ}="c",
\ar @{<-} "a";"b" <0pt>
\ar @{->} "b";"c" <0pt>
\endxy}\Ea
\  \to \
\Ba{c}\resizebox{11mm}{!}{
\xy
%
 (0,0)*{\circ}="a",
(8,0)*{\circ}="b",
\ar @{.} "a";"b" <0pt>
\endxy}\Ea
\]
\end{itemize}

It is elementary to check that $\pi_1$ and $\pi_2$ are both morphisms of dg Lie algebras.
The main Theorem is then an immediately consequence of the following two Propositions.

\sip

\subsubsection{\bf Proposition} {\em The map $\pi_2$ is a quasi-isomorphism.}
\begin{proof} Since $\pi_2$ is surjective, the statement of the proposition is equivalent to the fact that $\ker \pi_2=\fcGCdd$ is acyclic, which we shall show.
Let us consider a filtration of $\fcGCdd$ by the total number of vertices. The induced differential in the associated graded complex $gr\fcGCdd$
 just makes white vertices (if any) into black ones, $\circ \rightarrow \bu$.
By Maschke's Theorem, to prove the acyclicity of $gr\fcGCdd$ it is enough to prove the acyclicity
of its version $gr\fcGCdd^{marked}$ in which the generating graphs $\Ga$ have all their
edges and vertices totally ordered, but the type of vertices is not fixed. Then
there is an isomorphism of complexes
$$
gr\fcGCdd^{marked}=\prod_{\text{generators}\ \Ga}\left(\bigotimes_{v\in V(\Ga)} C_v\right)
$$
where $C_v$ is either a one-dimensional trivial complex (generated by a black vertex) in the
case when $v$ has at least one outgoing solid edge in ${\Ga}$, or $C_v$ is a two-dimensional
acyclic complex (generated by one black and one white vertex) if $v$ has no outgoing solid edges in $\Ga$. As any generating
graph $\Ga$  has least one vertex of the latter type since there cannot be directed cycles, the complex $gr\fcGCdd^{marked}$
has at least one acyclic tensor factor. The claim is proven
 (cf.\ Lemma 6.2.1 in \cite{Me0}).
\end{proof}

\sip

\subsubsection{\bf Proposition} {\em The map $\pi_1$ is a quasi-isomorphism.}
\begin{proof}
Note that the map $\pi_1$ is also surjective. Hence it is sufficient to check that 
\[
\hOGC_{d,d+1}^\bu := \ker \pi_1
\]
is acyclic. Furthermore note that $\hOGC_{d,d+1}^\bu\subset \hOGC_{d,d+1}$ is the subcomplex spanned by graphs with at least one black essential vertex.
Let us filter $\hOGC_{d,d+1}^\bu$ by the total number of essential vertices and consider the associated spectral sequence.
We claim that the associated graded with respect to this essential vertex filtration is acyclic
\begin{equation}\label{equ:tbs1}
H^\bu(gr_e \hOGC_{d,d+1}^\bu ) =0,
\end{equation}
from which the proposition immediately follows.
To show \eqref{equ:tbs1} we in turn consider on $gr_e \hOGC_{d,d+1}^\bu$ another filtration by the number of black vertices. Again it suffices to show that 
\begin{equation}\label{equ:tbs2}
H^\bu(gr_bgr_e \hOGC_{d,d+1}^\bu ) =0.
\end{equation}
To see this, note that the differential on the complex $gr_b gr_e \hOGC_{d,d+1}^\bu$ has only the following two terms, $\delta=\delta_1+\delta_2$ with 
\begin{itemize}
    \item $\delta_1$ is the natural inclusion of $\OGC_{d+1}^{0}$ into (the kernel of $\pi_1$ inside) $\fcGCdd$.
    \item $\delta_2$ creates one new inessential vertex from an edge between two essential vertices:
$$
\Ba{c}\resizebox{11mm}{!}{
\xy
%
 (0,0)*{\bullet}="a",
(8,0)*{\circ}="b",
\ar @{->} "a";"b" <0pt>
\endxy}\Ea
\  \to \
 \Ba{c}\resizebox{18mm}{!}{
\xy
%
 (0,0)*{\bullet}="a",
(8,0)*{\bullet}="b",
(16,0)*{\circ}="c",
\ar @{<-} "a";"b" <0pt>
\ar @{->} "b";"c" <0pt>
\endxy}\Ea
,
\quad
\quad
\quad 
\Ba{c}\resizebox{11mm}{!}{
\xy
%
 (0,0)*{\bullet}="a",
(8,0)*{\bullet}="b",
\ar @{->} "a";"b" <0pt>
\endxy}\Ea
\  \to \
 \Ba{c}\resizebox{18mm}{!}{
\xy
%
 (0,0)*{\bullet}="a",
(8,0)*{\bullet}="b",
(16,0)*{\bullet}="c",
\ar @{<-} "a";"b" <0pt>
\ar @{->} "b";"c" <0pt>
\endxy}\Ea.
$$
\end{itemize}
The map $\delta_1:  \OGC_{d+1}^{0} \to \fcGCdd$ is injective. 
Hence it is sufficient to check that 
\[
\mathrm{coker}(\delta_1)=: \OGC_{d,d+1}^{\bu\circ}
\]
is acyclic, i.e., 
\[
H(\OGC_{d,d+1}^{\bu\circ} ,\delta_2)=0.
\]
Note that $\OGC_{d,d+1}^{\bu\circ}$ can be identified with the subquotient of $\fcGCdd$ spanned by graphs with at least one white vertex and at least one black essential vertex.
Finally, the acyclicity of $(\OGC_{d,d+1}^{\bu\circ}, \delta_2)$ 
can be easily shown using the argument of \S 6.2.3 in \cite{Me0}. Indeed, one can assume 
without loss of generality that the essential vertices of graphs $\Ga$ generating 
$\OGC_{d,d+1}^{\bu\circ}$  are distinguished, say totally ordered. Every such graph 
$\Ga$ contains at least one black essential vertex and at least one white vertex which are connected 
by an edge, or are both neighbors of some (the same) inessential vertex, and we can assume without loss of generality that their labels 
are 1 and 2 respectively.  Considering a filtration of $\OGC_{d,d+1}^{\bu\circ}$ by 
the number of inessential vertices not between 1 and 2, we 
arrive at the associated graded complex which is the tensor product of a trivial complex 
and the complex $C_{12}$ which controls the types of all possible "edges" between 
vertices 1 and 2. One has
$$
C_{12}=\bigoplus_{k\geq 1} \odot^k C
$$
where $C$ is a 2-dimensional complex generated by the two possible connections,
$$
C=\text{span}\left\langle\Ba{c}\resizebox{11mm}{!}{
\xy
(0,2)*{^{1}},
(8,2)*{^{2}},
 (0,0)*{\bullet}="a",
(8,0)*{\circ}="b",
\ar @{->} "a";"b" <0pt>
\endxy}\Ea
\  ,\
 \Ba{c}\resizebox{18mm}{!}{
\xy
(0,2)*{^{1}},
(16,2)*{^{2}},
 (0,0)*{\bullet}="a",
(8,0)*{\bullet}="b",
(16,0)*{\circ}="c",
\ar @{<-} "a";"b" <0pt>
\ar @{->} "b";"c" <0pt>
\endxy}\Ea  \right\rangle
$$
with the differential sending the first element to the second. This complex is acyclic implying the acyclicity of $\OGC_{d,d+1}^{\bu\circ}$. The Proposition is proven.
\end{proof}

The proof of the Main Theorem is hence completed.

\sip

\def\cprime{$'$}

\end{document}